\documentclass[12pt]{article}
\usepackage{blindtext}
\usepackage[T1]{fontenc}
\usepackage[utf8]{inputenc}
\usepackage[greek, english]{babel}
\usepackage{hyperref}
\usepackage{inputenc}
\usepackage{amsfonts}
\usepackage{amssymb}
\usepackage{eucal}
\usepackage{amsxtra}
\usepackage{graphicx}
\usepackage{amsthm}
\usepackage{alphabeta}
\usepackage[margin=2cm]{geometry}
\usepackage{biblatex}
\theoremstyle{definition}

\newtheorem{theorem}{Theorem}[section]
\newtheorem{lem}{Lemma}[section]

\graphicspath{{./images/}}

\title{\textbf{A NOTE ON HARMONIC MAPS}}
\author{G. POLYCHROU}
\date{}
\begin{document}
\maketitle
\begin{abstract}
	\par In this note we will fill out the details from the recent work of Fotiadis and Daskaloyannis in \cite{2}, where the harmonic maps  described in \cite{5} are  written by the use of  Jacobi elliptic functions. We also prove that the harmonic maps described in \cite{6} can be written by the use of  Jacobi elliptic functions. 
\end{abstract}
\section{Introduction and statement of results}
	\par Let $(M,γ)$ and $(N,g)$ be Riemannian manifolds of dimension m and n, respectively. We use the following notations
	$(\gamma^{\alpha\beta}) = (\gamma_{\alpha\beta})^{-1}$,$\gamma = det(\gamma_{\alpha\beta})$ and
	$\Gamma^{\alpha}_{\beta\eta} = \frac{1}{2}\gamma^{\alpha\delta}(\gamma_{\beta\delta,\eta} + \gamma_{\eta\delta,\beta} - \gamma_{\beta\eta,\delta})$
	where $\Gamma^{\alpha}_{\beta\eta}$ are the Christoffel symbols of $Μ$. If $f:M \rightarrow N$ is a differentiable map with continuous derivative, we define its energy density with respect to the local coordinates $(x^{1},...,x^{m})$ on $M$ and $(f^{1},...,f^{n})$ on $N$ as 
		\begin{align}
		e(f)(x) = \frac{1}{2}\gamma^{\alpha\beta}(x)g_{ij}(f(x))\frac{\partial f^{i}(x)}{\partial x^{\alpha}}\frac{\partial f^{j}(x)}{\partial x^{\beta}}.
		\end{align}
	The energy of a map $f : M \rightarrow N$ is 
		\begin{align}
		E(f) = \int_{M} e(f)dM,
		\end{align}
	where $dM = \sqrt{\gamma}dx^{1}\wedge...\wedge dx^{m}$ in local coordinates is the volume form of $M$ and the Euler-Lagrange equations are
		\begin{align}
		\frac{1}{\sqrt{\gamma}}\frac{\partial}{\partial x^{\alpha}}(\sqrt{\gamma}\gamma^{\alpha\beta}\frac{\partial}{\partial x^{\beta}}f^{i}) + \gamma^{\alpha\beta}(x)\Gamma^{i}_{jk}(f(x))\frac{\partial}{\partial x^{\alpha}}f^{j}\frac{\partial}{\partial x^{\beta}}f^{k} = 0,
		\end{align}
	where  $α,β =1,...,m$ and $i,j,k = 1,...,n$. The proof can be found in  \cite{3}.
\par In this work we are interested in the case where $(M,g)$ and $(N,h)$ are Riemann surfaces. We introduce isothermal coordinates,
i.e. a  coordinate system $(x,y)$ on $M$  and $(R,S)$ on $N$ respectively such that
$$
g = e^{f(x,y)}(dx^{2} + dy^{2}) = e^{f(z,\bar{z})}dzd\bar{z} = e^{f(z,\bar{z})}|dz^{2}|,
$$
$$
h = e^{F(x,y)}(dR^{2} + dS^{2}) = e^{F(u,\bar{u})}dud\bar{u} = e^{F(u,\bar{u})}|du|^{2},
$$
where $z = x + iy$ and $u =R + iS$. In this case equations (3) turn into 
\begin{align}
\partial_{z\bar{z}}u + \partial_{u}F(u,\bar{u})\partial_{z}u\partial_{\bar{z}}u = 0,
\end{align}
where
$
\frac{\partial}{\partial z} = \frac{1}{2}(\frac{\partial}{\partial x} - i\frac{\partial}{\partial y})
$
and
$
\frac{\partial}{\partial \bar{z}} = \frac{1}{2}(\frac{\partial}{\partial x} + i\frac{\partial}{\partial y}).
$
\par In this note the example of harmonic maps presented in \cite{5} and \cite{6} are written in an analytic way by using the Jacobi elliptic functions. The main results are the following.
\begin{theorem}
		The  harmonic maps $u = (R,S)$ from the infinity strip  equipped with the conformal metric $dz = \frac{1}{\sin^{2}(y)}(dx^{2} + dy^{2})$ to itself of the form $R(x,y) = αx + h(y)$ and $S(x,y) = g(y)$ are
		\begin{align}
		S(x,y) = arccot ( wcs(αwy |\sin^{2}λ))
		\end{align}
	
		\begin{align}
		R(x,y) = αx + \frac{a^{2}y}{1 - w^{2}} - \frac{a^{2}}{αw(1 - w^{2})} Π(1 - \frac{1}{w^{2}};αwy|\sin^{2}λ)
		\end{align}
		where  $0 < λ < \frac{π}{2}$,
		\begin{align}	
		w = \sqrt{\frac{c^{2} + \sqrt{c^{4}- 4α^{2}b^{2}}}{2α^{2}}},
		\end{align}
	    $a = h'(\frac{\pi}{2})$ and $c^{2} = α^{2} + b^{2} + a^{4}.$
\end{theorem}
	
\begin{theorem}
	Let $u = (R,S)$ be a harmonic map from the infinity strip  equipped with the conformal metric $dz = \frac{1}{\sin^{2}(y)}(dx^{2} + dy^{2})$ to itself of the form $R(x,y) = αx$ and $S(x,y) = g(y)$, where $x$ is a real number , $0 \leq y \leq \frac{π}{2}$ and the boundary conditions $g(0) = 0$, $g(\frac{π}{2}) =\frac{π}{2}$ and $g'(\frac{π}{2}) =b$ hold true. Then,
\begin{align}
g(y) = arccot(cs(αy|1 - (\frac{b}{α})^{2})).	
\end{align}
	and the harmonic map $u$ can be extended for $0 \leq y \leq π$ with $g(π) = π$.
\end{theorem}
\section{Jacobi Elliptic functions}
\par The incomplete elliptic integral of the first kind $F$ is defined in \cite{1} as
	$$F(θ|m) = \int_{0}^{θ}\frac{dφ}{\sqrt{1 - m\sin^{2}(φ)}} = x.$$
	The real quarter period is $K = F(\frac{π}{2}|m)$ and $m$ is a parameter. Consider the complementary parameter $m'$ by $m + m' = 1$ and  define the imaginary quarter period $iK'$ as
	$$
	iK'(m) = iK' = i\int_{0}^{\frac{π}{2}}\frac{dφ}{\sqrt{1 - m'\sin^{2}(φ)}}.
	$$
	The incomplete elliptic integral of the third kind is defined as
	$$
	Π(n;θ|m) = \int_{0}^{θ}\frac{1}{1 - n\sin^{2}(φ)}\frac{dφ}{\sqrt{1 - m\sin^{2}(φ)}}
	= \int_{0}^{θ}\frac{du}{1 - n sn^{2}(u|m)}
	.
	$$
Let
\begin{align}
x = x(θ) = \int_{0}^{θ} \frac{dφ}{\sqrt{1 - m\sin^{2}(φ)}}.
\end{align}

\par Consider the three basic Jacobian elliptic functions:
\begin{align}
sn(x) = sn(x|m) = \sin(θ)
\end{align}
\begin{align}
cn(x) = cn(x|m) = \cos(θ)
\end{align}
\begin{align}
dn(x) = dn(x|m) = \sqrt{1 - m\sin^{2}(θ)} = \sqrt{1 - msn^{2}{x}}.
\end{align}
\par If $p,q,r$ are any three letters $s,c,d,n$ then  define 
\begin{align}
pq = \frac{pr}{qr} , pp = 1.
\end{align}
\par Then,
\begin{align}
	sn^{2}(z|m) + cn^{2}(z|m) = 1
\end{align}
\begin{align}
	sc(z + K|m) = -\frac{1}{\sqrt{m'}}cs(z|m)
\end{align}
\begin{align}	
	\frac{dcs(z)}{dz} = -\frac{dn(z)}{sn^{2}(z)},
\end{align}
where K is the quarter period and $m$ is a parameter and $m'$ is the complementary parameter corresponding $m$ and $m + m' = 1$. The proofs can be found in \cite{1} and in \cite{4}.
\par
 The inverses of the main three Jacobi elliptic functions can be represented as elliptic integrals.
 For example,we have
 \begin{align}
arcsn(x|m) = \int_{0}^{x}\frac{dt}{\sqrt{(1 - t^{2})(1 - mt^{2})}}.
\end{align}
	\section{The example of Shi, Tam and Wan \cite{5}}
	\par In this section we give a short introduction to the section 4 of \cite{5}. Consider the infinity strip $\mathbb{H}^{2} = \{(x,y) \in \mathbb{R}^{2} | y \in (0,π) \}$ equipped with the conformal metric $dz = \frac{1}{\sin^{2}(y)}(dx^{2} + dy^{2})$. This model of the hyperbolic space is related to one-soliton solutions of the  Sinh-Gordon equation as it has been proved in \cite{2}. The Euler-Lagrange equations of a harmonic map $u(x,y) = (R,S)$ from the infinite strip to itself are as follows:
	\begin{align}	
	\bigtriangleup_{0}R - 2\cot(S)<\bigtriangledown_{0}R,\bigtriangledown_{0}S>=0
	\end{align}
	\begin{align}	
	\bigtriangleup_{0}S + \cot(S)(<\bigtriangledown_{0}R,\bigtriangledown_{0}R> - <\bigtriangledown_{0}S,\bigtriangledown_{0}S>)=0,
	\end{align}
	where $\bigtriangleup_{0} = \frac{\partial^{2}}{\partial x^{2}} +  \frac{\partial^{2}}{\partial y^{2}}$ is the euclidean Laplacian and $\bigtriangledown_{0}$ the euclidean gradient. In \cite{5} they consider the boundary values 
	\begin{align}
	R(x,0) = αx, R(x,π) = αx+β, \; β > 0
	\end{align}
	\begin{align}
	S(x,0)=0, S(x,π)=π
	\end{align}
	 and $y$ takes values on the interval $(0,π)$  .  In this case the Euler-Lagrange equations become:
	\begin{align}	
	h'' - 2\cot(g)g'h' = 0 
	\end{align}
	\begin{align}	
	g'' + \cot(g)(α^{2} + (h')^{2} - (g')^{2}) = 0,
	\end{align}
	with boundary values $h(0) = 0, h(π) = β, g(0) = 0, g(π) = π$ and $g$ takes values on the interval $(0,π)$. Since we  have solutions with $h' \geq 0$, we get that
	\begin{align}	
	\frac{\partial R}{\partial y} = h' = a^{2}\sin^{2}(g),
	\end{align}
	\begin{align}	
	\frac{\partial S}{\partial y} = (g')^{2} = α^{2} + (b^{2} + a^{4} - α^{2})\sin^{2}(g) - a^{4}\sin^{4}(g),
	\end{align}
	where 
	\begin{align}
	a = h'(\frac{π}{2})
	\end{align} 
	and 
	\begin{align}
	b = g'(\frac{π}{2}) 
	\end{align} 
	are constants to be chosen such that $u = (R,S)$ satisfies the boundary values. In \cite{5} they study harmonic maps of the form $u = (R , S)$ and $R(x,y) = αx + h(y)$ and $S(x,y) = g(y)$  ,both of $h$ and $g$ are defined on the interval $[0,\frac{π}{2}]$. The functions $g$ and $h$ can be extended in $[0,π]$ using the following formulas:
	$$
	h(y) = β - h(π - y), g(y) = π - g(π - y).
	$$
	\par In order to simplify equations (24) and (25) we let $z = \cot(g)$. By the fact that \\$\sin^{2}(g) = \frac{1}{1 + \cot^{2}(g)} = \frac{1}{1 + z^{2}}$ and $g' = z'\sin^{2}(g) = \frac{z'}{1 + z^{2}}$, we conclude that
	$$
	(\frac{z'}{1 + z^{2}})^{2} = α^{2} + (b^{2} + a^{4} - α^{2})\frac{1}{1 + z^{2}} - a^{4}(\frac{1}{1 + z^{2}})^{2}.
	$$
	Then $(z')^{2} = α^{2}z^{4} + (α^{2} + b^{2} + a^{4})z^{2} +b^{2}$.  So $z(y)$ is defined by
	\begin{gather}
	\int_{0}^{z(y)} \frac{dz}{\sqrt{α^{2}z^{4} + (α^{2} + b^{2} + a^{4})z^{2} +b^{2}}} = \frac{π}{2} - y.
	\end{gather}
	\par In \cite{5}, Shi, Tam and Wan have proved the two following lemmas.
	\begin{lem}
		For any $a \geq 0$, there is a unique constant $b_{a} > 0$ such that 
		\begin{align}	
		\int_{0}^{\infty}\frac{dz}{\sqrt{α^{2}z^{4} + (α^{2} + b_{a}^{2} + a^{4})z^{2} +b_{a}^{2}}} = \frac{π}{2}.
		\end{align}
		Moreover, $b_{a} \leq \max\{2α,α^{-1}\}$ and $b_{a}$ depends continuously on a.
	\end{lem}
	\begin{lem}
		 Given $α > 0$ and $β \geq 0$ there is $a = a({α,β}) \geq 0$ such that
		\begin{align}	
		a^{2}\int_{0}^{\infty}\frac{dz}{(1+z^{2})\sqrt{α^{2}z^{4} + (α^{2} + b_{a}^{2} + a^{4})z^{2} +b_{a}^{2}}} = \frac{β}{2}.
		\end{align}
	\end{lem}
	\section{The solution by the use of  elliptic functions}
	\par  Using elliptic functions we study the relations between the constants $α,a,b,c$ and $β$. Moreover using elliptic functions we prove relations, which describe the quarter period $K(\sin^{2}λ)$ and the extension of the solutions from the interval $[0,\frac{π}{2}]$ to $[0,π]$ follows easier by using properties of the elliptic functions. 
	\par In this section consider the boundary conditions $R(x,0) = αx$, $R(x,π) = αx + β$ where $β > 0$  which implies that $a > 0$, $S(x,0) = 0$ and $S(x,π) = π$. The main result is the following.
\begin{theorem}
		The  harmonic maps $u = (R,S)$ from the infinity strip  equipped with the conformal metric $dz = \frac{1}{\sin^{2}(y)}(dx^{2} + dy^{2})$ to itself of the form $R(x,y) = αx + h(y)$ and $S(x,y) = g(y)$ are
		\begin{align}
		S(x,y) = arccot ( wcs(αwy |\sin^{2}λ))
		\end{align}
		\begin{align}
		R(x,y) = αx + \frac{a^{2}y}{1 - w^{2}} - \frac{a^{2}}{αw(1 - w^{2})} Π(1 - \frac{1}{w^{2}};αwy|\sin^{2}λ)
		\end{align}
		where  $0 < λ < \frac{π}{2}$,
		\begin{align}	
		w = \sqrt{\frac{c^{2} + \sqrt{c^{4}- 4α^{2}b^{2}}}{2α^{2}}},
		\end{align}
		 $a = h'(\frac{\pi}{2})$ and $c^{2} = α^{2} + b^{2} + a^{4}.$
\end{theorem}
	\begin{proof}
		\par Consider  $z(y)$ as in equation (28). It holds that,
		$$
		α^{2}z^{4} + c^{2}z^{2} + b^{2} = α^{2}(z^{2} + (w\cos λ)^{2})(z^{2} + w^{2}),
		$$
		where $$w = \sqrt{\frac{c^{2} + \sqrt{c^{4}- 4α^{2}b^{2}}}{2α^{2}}}.$$ 
		So it holds that
$$
		Ι = \frac{π}{2} - y = \int_{0}^{z(y)}\frac{dz}{\sqrt{α^{2}z^{4} + c^{2}z^{2} + b^{2}}} = \frac{1}{α}\int_{0}^{z(y)}\frac{dz}{\sqrt{(z^{2} + (w\cos λ)^{2}))(z^{2} + w^{2})}}.
$$
		Now, let $u^{2} = \frac{z^{2}}{z^{2} + (w\cos λ)^{2}}$.  We obtain the following:
		$$
		\frac{π}{2} - y = \frac{1}{αw}\int_{0}^{\sqrt{\frac{z^{2}}{z^{2} + (w\cos λ)^{2}}}} \frac{du}{\sqrt{(1 - u^{2})(1 - \sin^{2}λ u^{2})}}.
		$$
		By the definition of the elliptic sine function $sn$ we get that
		$$
		sn(αw\frac{π}{2} - αwy|\sin^{2}λ) = \sqrt{\frac{z^{2}(y)}{z^{2}(y) + (w\cos λ)^{2}}},
		$$
		and 
		$$
		z^{2}(y) = \frac{(w\cos λ)^{2}sn^{2}(αw\frac{π}{2} - αwy|\sin^{2}λ)}{1 - sn^{2}(αw\frac{π}{2} - αwy|\sin^{2}λ)}.
		$$
		Thus,
		\begin{align}
		z(y) = w\cos λ sc(αw\frac{π}{2} - αwy|\sin^{2}λ).
		\end{align}
		By Lemma 3.1, there is a unique constant $b_{α}$ such that when $y$ tends to zero then $z(y)$ tends to infinity. 
			So by the property $sc(u) = \frac{sn(u)}{cn(u)}$ we get that: 
			$$
			cn(αw\frac{π}{2}|\sin^{2}λ) = 0.
			$$
			Using the fact that $cn(K|\sin^{2}λ) = 0$, where $K(\sin^{2}λ)$ is the quarter period we obtain that ,
			\begin{align}
			K(\sin^{2}λ) = αw\frac{π}{2}.
			\end{align} 
		By (35),  equation (34) turns into
		$$
		z(y) = w\cos λ sc(αw\frac{π}{2} - αwy|\sin^{2}λ) = w\cos λ sc(Κ - αwy|\sin^{2}λ), 
		$$
		thus,
		\begin{align}
		z(y) = wcs(αwy|\sin^{2}λ).
		\end{align}
		So, using the identity (16)
		we see that,
		\begin{align}
		z'(y) = -αw^{2}\frac{dn(αwy|\sin^{2}λ)}{sn^{2}(αwy|\sin^{2}λ)}.
		\end{align}
		 Moreover,
		$$
		z(y) = \cot(g(y)) \Rightarrow g(y) = arccotz(y) \Rightarrow S(x,y)= arccotz(y),
		$$
		thus
		\begin{align}	
		S(x,y)= arccot(wcs(αwy|\sin^{2}λ))).
		\end{align}

		By  (38) and using equation (37) we conclude that
		$$
		\frac{\partial S}{\partial y} = -\frac{z'(y)}{1 + z^{2}(y)} = \frac{αw^{2}dn(αwy|\sin^{2}λ)}{sn^{2}(αwy|\sin^{2}λ)	+ w^{2}cn^{2}(αwy|\sin^{2}λ)},
		$$ 
		and by using (14) we have that
		\begin{align}	
		\frac{\partial S}{\partial y} =  \frac{αw^{2}dn(αwy|\sin^{2}λ)}{w^{2} + (1 - w^{2}) sn^{2}(αwy|\sin^{2}λ)}.
		\end{align}
		We also have,
		$$
		\frac{\partial R}{\partial y}  = a^{2}sin^{2}(S) = a^{2}\sin^{2}(arccot(wcs(αwy|\sin^{2}λ))) = a^{2}\frac{1}{1 + w^{2}cs^{2}(αwy|\sin^{2}λ)}.
		$$
		Thus,
		\begin{align}	
		\frac{\partial R}{\partial y} = a^{2}\frac{1}{1 + w^{2}cs^{2}(αwy|\sin^{2}λ)} = a^{2}\frac{sn^{2}(αwy|\sin^{2}λ)}{w^{2} + (1 - w^{2}) sn^{2}(αwy|\sin^{2}λ)}
		\end{align}	
and so we conclude that
\begin{align}
	h'(y) = a^{2}\frac{sn^{2}(αwy|\sin^{2}λ)}{w^{2} + (1 - w^{2}) sn^{2}(αwy|\sin^{2}λ)}
\end{align}
$$
\Rightarrow h(y) = \int_{0}^{y}a^{2}\frac{sn^{2}(αwu|\sin^{2}λ)}
{w^{2} - (w^{2} - 1) sn^{2}(αwu|\sin^{2}λ)}du 
=\frac{a^{2}}{1 - w^{2}}\int_{0}^{y}
\frac{-(w^{2} - 1)sn^{2}(αwu|\sin^{2}λ)}{w^{2} - (w^{2} - 1) sn^{2}(αwu|\sin^{2}λ)}du 
$$
$$
= \frac{a^{2}}{1 - w^{2}}\int_{0}^{y}
\frac{w^{2} - (w^{2} - 1)sn^{2}(αwu|\sin^{2}λ) - w^{2}}
{w^{2} - (w^{2} - 1) sn^{2}(αwu|\sin^{2}λ)}du
$$
$$
= \frac{a^{2}}{1 - w^{2}}\int_{0}^{y}du
- \frac{a^{2}w^{2}}{1 - w^{2}}\int_{0}^{y}
\frac{1}{w^{2} - (w^{2} - 1) sn^{2}(αwu|\sin^{2}λ)}du 
$$
$$
= \frac{a^{2}}{1 - w^{2}}\int_{0}^{y}du
- \frac{a^{2}}{1 - w^{2}}\int_{0}^{y}
\frac{1}{1 - (1 - \frac{1}{w^{2}}) sn^{2}(αwu|\sin^{2}λ)}du
$$
\begin{align}
	\Rightarrow
	h(y) = \frac{a^{2}y}{1 - w^{2}} - \frac{a^{2}}{αw(1 - w^{2})} Π(1 - \frac{1}{w^{2}};αwy|\sin^{2}λ)
\end{align}
and so we figure out that
\begin{align}
R(x,y) = αx + \frac{a^{2}y}{1 - w^{2}} - \frac{a^{2}}{αw(1 - w^{2})} Π(1 - \frac{1}{w^{2}};αwy|\sin^{2}λ).
\end{align}
\end{proof}
		\begin{lem}
		The following equation hold true
			$$
			\frac{β}{2} = a^{2}\int_{0}^{\infty}\frac{dz}{(1+z^{2})\sqrt{α^{2}z^{4} + (α^{2} + b^{2} + a^{4})z^{2} +b^{2}}} = h(\frac{\pi}{2})
			$$ 
			where $w$ and $c$ as before.
		\end{lem}
		\begin{proof}
		$$
		I = a^{2}\int_{0}^{\infty}\frac{dz}{(1+z^{2})\sqrt{α^{2}z^{4} + (α^{2} + b^{2} + a^{4})z^{2} +b^{2}}}
		= \frac{a^{2}}{α}\int_{0}^{\infty}\frac{dz}{(1+z^{2})\sqrt{(z^{2} + (w\cos λ)^{2})(z^{2} + w^{2})}}.
		$$
		We apply the change of variables $z = \cot φ$, so
		$$
		I = \frac{a^{2}w}{α}\int_{0}^{\frac{\pi}{2}}\frac{dφ}
		{\sin^{2}φ(1 + w^{2}\cot^{2}φ)\sqrt{(w^{2}\cot^{2}φ + w^{2}\cos^{2}λ)(w^{2}\cot^{2}φ + w^{2})}}.
		$$
		$$
		\Rightarrow
		I = \frac{a^{2}}{αw}\int_{0}^{\frac{\pi}{2}}\frac{\sin^{2}φ}{w^{2} - (w^{2} - 1)\sin^{2}φ}
		\frac{dφ}{\sqrt{1 - \sin^{2}λ\sin^{2}φ}}.
		$$
		By letting $sn(αwy|\sin^{2}λ) = \sin φ$ then $dφ = αw \sqrt{1 - \sin^{2}λ\sin^{2}φ}dy$ we get that
		$$
		\Rightarrow I =
		a^{2} \int_{0}^{\frac{\pi}{2}} \frac{sn^{2}(αwy|\sin^{2}λ)}{w^{2} - (w^{2} - 1)sn^{2}(αwy|\sin^{2}λ)}dy.
		$$
		and using equation (41) we obtain the requested.
		\end{proof}

		\par Furthermore, one can easily check that the symmetry of the solutions described in \cite{5} follows by using the following properties of the Jacobi elliptic functions
		$$
		dn(αw(π - y)) = dn(2K - αwy) =  dn(αwy),
		$$ 
		since $dn(z)$ is an even function and 
		$$
		sn(αw(π - y)) = sn(2K - αwy) =   sn(αwy).
		$$
		Note that $sn(z)$ is an odd function thus we obtain that $g'(π - y) = g'(y)$ and $h'(π - y) = h'(y)$ and since $g(\frac{π}{2}) = \frac{π}{2}$ and $h(\frac{π}{2}) = \frac{β}{2}$ from the definition of $h$ we see that the followings hold true:
		$$
		h(t)=β-h(π-t) , g(t)=π-g(π-t).
		$$
	\section{The example of Wang \cite{6}}
	\par Now we are going to solve the boundary value problem when $β = 0$. So the boundary conditions for the function $R$ are $R(x,0) = 0$ and $R(x,π) = 0$, while $S$ as in Section 3. So $h(y) = 0$ and $g(0) = 0$, $g(\frac{π}{2}) =\frac{π}{2}$, $g(π) = π$ and $g'(\frac{π}{2}) =b$. Then the solutions  are $R(x,y) = αx$ and $S(x,y) = g(y)$. If $β = 0$, then we have that $a = 0$. So this case differs from the one described in Section 4. In this case the Euler-Lagrange equations turn into a non-linear ordinary differential equation:
	\begin{align}
g'' + \cot(g)(α^{2} - (g')^{2}) = 0 ,
	\end{align}
	\par see \cite{3}. The main result of this section is the following.
	\begin{theorem}
		The function $g(y)$ described above is  
	\begin{align}
	g(y) = arccot(cs(αy|1 - (\frac{b}{α})^{2})).
	\end{align}
	\end{theorem}
\begin{proof}
	Equation (44) can be written as
	$$
	\frac{2g'g''}{(g')^{2} - α^{2}} = 2g'\cot(g).
	$$
	By integrating both sides we get that
	$$
	\ln(|(g')^{2} - α^{2}|) = 2\ln(\sin(g)) + \ln(c).
	$$
	Using of the boundary values we have that
	$$
	(g')^{2} -  α^{2} = c\sin^{2}(g),
	$$
	where $c = b^{2} - α^{2}$. So 
	$$
	g' = \sqrt{α^{2} - (α^{2} - b^{2})\sin^{2}(g)} \Leftrightarrow \frac{g'(y)}{α} = \sqrt{1 - (1 - \frac{b^{2}}{α^{2}})\sin^{2}(g)}.
	$$
	Thus,
	$$
	αdy = \frac{dg}{\sqrt{1 - (1 - \frac{b^{2}}{α^{2}})\sin^{2}(g)}} \Rightarrow \int_{0}^{y}αdy = \int_{0}^{g}\frac{du}{\sqrt{1 - (1 - \frac{b^{2}}{α^{2}})\sin^{2}(u)}}.
	$$
	$$
	\Rightarrow αy = \int_{0}^{g}\frac{du}{\sqrt{1 - (1 - \frac{b^{2}}{α^{2}})\sin^{2}(u)}}.
	$$
	We conclude that
	$$
	\sin(g(y)) = sn(αy|1 - \frac{b^{2}}{α^{2}}).
	$$
	So 
	\begin{align}
		g(y) = \arcsin(sn(αy|1 - (\frac{b}{α})^{2}))
	\end{align}
	and by the fact that $\arcsin(sn(x)) = arccot(cs(x))$ we obtain the requested
	\end{proof}
	\par Note that the case $β = 0$ can be solved by using the fact that $a = 0$ and making the substitution in Theorem 4.1. One can check that when $a = 0$ then $w = 1$ and $\sin^{2}λ =  1 - \frac{b^{2}}{α^{2}}$. 
	
	E-mail address: G. Polychrou, ipolychr@math.auth.gr
	Current address: Department of Mathematics, Aristotle University of Thessaloniki, Thessaloniki 54124, Greece
\end{document}